\documentclass[conference]{IEEEtran}
\IEEEoverridecommandlockouts
\pdfoutput=1 
\usepackage{cite}
\usepackage{amsmath,amssymb,amsfonts}
\usepackage{textcomp}
\usepackage{xcolor}

\usepackage{bm}
\usepackage{amsthm}
\usepackage{mathrsfs}
\usepackage{algorithm}
\usepackage{algorithmicx}
\usepackage{algpseudocode}
\usepackage{mathabx}
\usepackage{enumerate}

\usepackage{graphicx}
\usepackage[hang,small]{caption}
\usepackage[subrefformat=parens]{subcaption}

\algdef{SE}[DOWHILE]{Do}{doWhile}{\algorithmicdo}[1]{\algorithmicwhile\ #1}

\newcommand{\ip}[2]{\left\langle #1 , #2 \right\rangle}
\newcommand{\norm}[1]{\left\| #1 \right\|}
\newcommand{\minimize}[0]{\operatornamewithlimits{minimize}}

\newcommand{\ri}[0]{\operatorname{ri}}

\newcommand{\dom}{\operatorname{dom}}
\newcommand{\ran}{\operatorname{ran}}

\newcommand{\rank}{\operatorname{rank}}
\newcommand{\nullsp}{\operatorname{null}}

\newcommand{\lrangle}[1]{\left\langle #1 \right\rangle}
\newcommand{\argmin}{\operatornamewithlimits{argmin}}
\newcommand{\Id}{\operatorname{Id}}
\newcommand{\prox}{\operatorname{Prox}}
\newcommand{\fix}{\operatorname{Fix}}
\newcommand{\setLO}[2]{\mathcal{B}\left( #1 , #2 \right)}
\newcommand{\zeroMatrix}{\mathrm{O}}
\newcommand{\identityMatrix}{\mathrm{I}}

\newcommand{\lev}{\operatorname{lev}}
\newcommand{\setN}{\mathbb{N}}
\newcommand{\setR}{\mathbb{R}}
\newcommand{\setPR}{\mathbb{R}_{++}}
\newcommand{\setNNR}{\mathbb{R}_{+}}
\newcommand{\transT}{\mathsf{T}}
\newcommand{\vecrize}{\operatorname{vec}}
\newcommand{\ind}{\Upsilon}

\newcommand{\spX}{\mathcal{X}}
\newcommand{\spY}{\mathcal{Y}}
\newcommand{\spZ}{\mathcal{Z}}
\newcommand{\sptildeZ}{\tilde{\mathcal{Z}}}
\newcommand{\spfrakZ}{\mathfrak{Z}}
\newcommand{\spH}{\mathcal{H}}
\newcommand{\spK}{\mathcal{K}}
\newcommand{\TCC}{T_\mathrm{cLiGME}}
\newcommand{\opP}{\mathfrak{P}}
\newcommand{\opL}{\mathfrak{L}}

\newcommand{\opLC}{\opL_c}
\newcommand{\opC}{\mathfrak{C}}
\newcommand{\tildeL}{\tilde{\opL}}
\newcommand{\xstar}{x^\star}

\newcommand{\setS}{\mathcal{S}}

\newcommand{\bfPsi}{\mathbf{\Psi}}
\newcommand{\dmatv}{D_{\mathrm{V}}}
\newcommand{\dmath}{D_{\mathrm{H}}}

\newcommand{\RomanOne}{\mathrm{I}}
\newcommand{\RomanTwo}{\mathrm{I\hspace{-.1em}I}}
\newcommand{\RomanThree}{\mathrm{I\hspace{-.1em}I\hspace{-.1em}I}}
\newcommand{\RomanFour}{\mathrm{I\hspace{-.1em}V}}
\newcommand{\RomanFive}{\mathrm{V}}
\newcommand{\RomanSix}{\mathrm{V\hspace{-.1em}I}}
\newcommand{\RomanSeven}{\mathrm{V\hspace{-1pt}I\hspace{-1pt}I}}
\newcommand{\RomanEight}{\mathrm{V\hspace{-1pt}I\hspace{-1pt}I\hspace{-1pt}I}}

\newtheorem{theorem}{Theorem}[]

\newtheorem{lemma}{Lemma}[]
\newtheorem{proposition}{Proposition}[]

\newtheorem{problem}{Problem}[]
\newtheorem{remark}{Remark}[]

\def\BibTeX{{\rm B\kern-.05em{\sc i\kern-.025em b}\kern-.08em
    T\kern-.1667em\lower.7ex\hbox{E}\kern-.125emX}}
\begin{document}

\title{A Convexly Constrained LiGME Model and\\
Its Proximal Splitting Algorithm
\thanks{This work was supported in part by JSPS Grants-in-Aid (19H04134).}
}

\author{\IEEEauthorblockN{Wataru Yata, Masao Yamagishi, and Isao Yamada}
\IEEEauthorblockA{\textit{Dept. of Information and Communications Engineering} \\
\textit{Tokyo Institute of Technology}, Tokyo, Japan\\
Email:\{yata,myamagi,isao\}@sp.ict.e.titech.ac.jp}
}

\maketitle
\begin{abstract}
For the sparsity-rank-aware least squares estimations,
the LiGME 
(Linearly involved Generalized Moreau Enhanced) model
was established recently
in [Abe, Yamagishi, Yamada, 2020]
to use certain nonconvex enhancements of
linearly involved
convex regularizers
without losing their overall convexities.
In this paper, 
for further advancement
of the LiGME model by incorporating
multiple a priori knowledge as hard convex constraints, 
we newly propose a convexly constrained LiGME (cLiGME) model.
The cLiGME model 
can utilize multiple convex constraints
while preserving benefits
achieved by the LiGME model.
We also present a proximal splitting type algorithm for
the proposed cLiGME model.
Numerical experiments demonstrate the efficacy of the proposed model 
and the proposed optimization algorithm
in a scenario of signal processing application.
\end{abstract}


\section{Introduction}
Many tasks in data science have been formulated as
an estimation of an unknown vector $\xstar\in \spX$
from the observed vector $y\in\spY$ which follows 
the linear regression model:
\begin{align}
  \label{eq: lrm}
  y = A\xstar + \varepsilon,
\end{align}
where
$(\spX,\ip{\cdot}{\cdot}_\spX,\norm{\cdot}_\spX)$
and
$(\spY, \ip{\cdot}{\cdot}_\spY, \norm{\cdot}_\spY)$
are finite dimensional real Hilbert spaces,
$A:\spX\to\spY$ is a known linear operator
and $\varepsilon$ is an unknown noise vector.
To tackle such estimations, we commonly formulate the following 
regularized 
least-squares problem:
\begin{align}
  \label{eq: regularized least square}
  \minimize_{x\in \spX}
  J_{\Psi\circ \opL}(x):=\frac{1}{2}\norm{Ax-y}_\spY^2
  +\mu \Psi\circ\opL(x),
\end{align}
where a function $\Psi:\spZ\to (-\infty,\infty]$
is designed over a finite dimensional real Hilbert space 
$(\spZ,\ip{\cdot}{\cdot}_\spZ,\norm{\cdot}_\spZ)$,
$\opL:\spX\to\spZ$ is a linear operator
and
$\mu\geq 0$ is a regularization parameter.
In the problem \eqref{eq: regularized least square},
$\frac{1}{2}\norm{A\cdot-y}_\spY^2$ is called the 
least-squares term and
$\Psi\circ\opL$ is called a regularizer (or penalty)
which is designed based on a priori knowledge regarding $\xstar\in\spX$.
The tuple $(\Psi, \opL)$ is designed not only 
depending on each application but also
based on mathematical tractability
of the optimization model \eqref{eq: regularized least square}.
\par

To design $(\Psi, \opL)$ in \eqref{eq: regularized least square},
one of the most crucial properties
to be cared in modern signal processing is the sparsity 
of vectors or the low-rankness of matrices.
The sparsity and
the low-rankness
have been used for many scenarios
in data science 
such as
compressed sensing
\cite{Stable signal recovery from incomplete and 
inaccurate measurements.,
Compressed sensing.}
and related sparsity
aware applications
\cite{Sparse and Redundant Representations.,
  Machine Learning: A Bayesian and 
  Optimization Perspective.}.
A naive measure of the sparsity
is the $l_0$-pseudonorm $\norm{\cdot}_0$
which stands for the number of nonzero components of a
given vector.
However, the problem 
\eqref{eq: regularized least square}
with $\Psi := \norm{\cdot}_0$ is known to be
NP-hard
\cite{Sparse approximate solutions to linear systems.,
NP-hardness of l0 minimization problems: revision and extension to the non-negative setting}.
To avoid this burden,
the $l_1$-norm $\norm{\cdot}_1$, which is the 
convex envelope of $\norm{\cdot}_0$
in the vicinity of zero vector,
has been used as $\Psi$ in many applications,
e.g., Lasso 
\cite{Regression shrinkage and selection via the lasso.}
with $(\Psi, \opL) := (\norm{\cdot}_1,\Id)$,
the convex total variation (TV) \cite{Nonlinear total variation based noise removal algorithms}
with $(\Psi,\opL):=(\norm{\cdot}_ 1,D)$
and the wavelet-based regularization 
\cite{An Iterative Thresholding Algorithm for Linear Inverse Problems with a Sparsity Constrains,
Sparse Image and Signal Processing}
with $(\Psi,\opL) := (\norm{\cdot}_1, W)$,
where 
$\Id$ is the identity operator,
$D$ is the first order differential operator 
and 
$W$ is a wavelet transform matrix.\par 

The convexity of $J_{\Psi\circ \opL}$
is a key property to obtain a global minimizer
of $J_{\Psi\circ \opL}$ even
by the state-of-the-art optimization techniques.
If $\Psi$ is a convex function,
$J_{\Psi \circ \opL}$ remains to be convex
because the sum of convex functions
remains convex.
However, 
the convexity of $\Psi$ is just a
sufficient condition 
for the overall convexity of $J_{\Psi\circ\opL}$.
In fact, we can design a nonconvex function $\Psi$
for $J_{\Psi\circ\opL}$ to be convex.
Such  penalties are called \textit{the
convexity-preserving nonconvex penalties}
\cite{Visual Reconstruction.,
Estimation of binary images by minimizing convex criteria.,
Markovian reconstruction using a GNC approach.,
Energy minimization methods.,
GMC,
Nearly unbiased variable selection under minimax concave
penalty.,
Stable principal component pursuit via convex analysis.
}.
To design the convexity-preserving nonconvex penalty,
the Linearly involved
Generalized Moreau Enhanced (LiGME) model 
\begin{align}
  \label{eq: LiGME}
  \minimize_{x\in\spX}
  J_{\Psi_B\circ \opL}(x):=\frac{1}{2}
  \norm{Ax-y}_{\spY}^2 + \mu \Psi_B\circ\opL(x)
\end{align}
was introduced in \cite{LiGME}
for a proper lower semicontinuous convex
function\footnote{
  A convex function $\Psi:\spZ\to(-\infty,\infty]$
  is said to be  $\Psi \in \Gamma_0(\spZ)$
  if $\Psi$ is proper 
  (i.e., $\dom(\Psi):=\{z\in\spZ|
  \Psi(z) <\infty\}
  \neq \varnothing$) and
  lower semicontinuous 
  (i.e., 
  $(\forall a\in \setR)\ 
  \lev_{\leq a} \Psi :=
  \{
    z\in \spZ |
    \Psi(z)\leq a
  \}
  $
  is closed).\par
  For $\Psi\in \Gamma_0(\spZ)$, 
  the proximity operator of $\Psi$ 
  is defined by 
  \begin{align*}
    \prox_\Psi :\spZ\to\spZ:x\mapsto \argmin_{y\in \spZ}
    \left[\Psi(y) + \frac{1}{2}\norm{x-y}_\spZ^2\right].
  \end{align*}
}
$\Psi\in \Gamma_0(\spZ)$,
satisfying coercivity 
(i.e., 
$\lim_{\norm{z}_\spZ\to\infty}\Psi(z)= \infty$) and
$\dom(\Psi) = \spZ$,
where $\Psi_B$ is the generalized Moreau
enhanced (GME) penalty
\begin{align}
  \Psi_B(\cdot) := \Psi(\cdot) - \min_{v\in \spZ}\left[ \Psi(v) 
  +\frac{1}{2}\norm{B\left(\cdot - v\right)}_{\sptildeZ}^2\right]
\end{align}
with a finite dimensional real
Hilbert space $\sptildeZ$ and a linear operator 
$B:\spZ\to\sptildeZ$.
Indeed, if 
$A^*A -\mu \opL^*B^*B\opL\succeq \zeroMatrix_\spX$,
the overall convexity of $J_{\Psi_B\circ \opL}$ is
achieved.
Compared with its special case 
$(\Psi,\opL) = (\norm{\cdot}_1,\Id)$
\cite{GMC,
Nearly unbiased variable selection under minimax concave
penalty.
}
and variants 
\cite{Stable principal component pursuit via convex analysis.},
the LiGME model can admit multiple 
linearly involved nonconvex
penalties,
by formulating it in a product space,
without losing the overall
convexity of $J_{\Psi_B\circ \opL}$ and thus is applicable to
broader range of scenarios in the sparsity-rank-aware 
signal processing.
\par 

In contrast, as seen in the
set theoretic estimation 
\cite{Parallel Optimization: Theory Algorithm and Optimization.,
Foundation of set theoretic estimation.,
Image Recovery: Theory and Application.,
Image Restoration by the Method of Convex Projections: Part 1-Theory.},
if we know a priori that our target vector $\xstar$ 
in \eqref{eq: lrm}
satisfies $\opC_i\xstar\in C_i\ (i\in I)$,
where $I$
is a finite index set,
$C_i\ (i\in I)$ is a closed convex subset of 
a finite dimensional real Hilbert space $\spfrakZ_i$ and
$\opC_i:\spX\to \spfrakZ_i$ is a linear operator,
we can expect to improve the estimation achieved
via the model \eqref{eq: LiGME}
by incorporating
multiple convex constraints
as 
\begin{align}
  \label{eq: LiGME with convex constraints}
  \hspace*{-3.5mm}\minimize_{( i \in I)\
  \opC_i x\in C_i\subset \spfrakZ_i}
  J_{\Psi_B\circ \opL}(x):=\frac{1}{2}\norm{Ax-y}_\spY^2
  +\mu \Psi_B\circ\opL(x).
\end{align}
Note that,
as we will see in Lemma \ref{def: cLiGME},
the model \eqref{eq: LiGME} itself
does not cover
the model \eqref{eq: LiGME with convex constraints}
because the condition $\dom(\Psi) = \spZ$ does not 
allow to
use, 
in the model \eqref{eq: LiGME},
the indicator function
\begin{align*}
  \iota_{C_i}(x):= \begin{cases}
    0 & (x\in C_i)\\
    \infty & (otherwise).
  \end{cases}
\end{align*}
\par

In this paper, 
we call the problem
\eqref{eq: LiGME with convex constraints}
a convexly constrained
Linearly involved Generalized Moreau Enhanced
(cLiGME) model.
This problem can be reformulated as 
simpler models
\eqref{eq: LiGME with C} and
\eqref{cLiGME}.
The overall convexity of
$J_{\bfPsi_{B_c}\circ\opLC}$
is guaranteed, 
by Proposition \ref{prop: overall convexity},
under a certain condition regarding a tunable matrix $B$.
We also present a way to design such a matrix $B$ satisfying the 
overall convexity condition
(see Proposition \ref{prop: how to choose B} 
and Proposition \ref{prop: multiple penalties}).
Moreover,
as an extension of \cite[Theorem 1]{LiGME},
we present a novel 
proximal
splitting type algorithm 
of guaranteed convergence to a global minimizer 
of $J_{\bfPsi_{B_c}\circ \opL_c}$
under its overall convexity condition 
(see Theorem \ref{thm: TCC}).
Finally, to demonstrate the efficacy of the proposed model
and algorithm, 
we present
in Section \ref{sec: numerical experiment} numerical experiments in a
simple sparsity-aware signal processing scenario.

\section{Preliminaries}
Let $\setN$, $\setR$, $\setNNR$ and $\setPR$
be the set of nonnegative integers, real numbers,
nonnegative real numbers and positive real numbers, respectively.
Let
$(\spH,
\ip{\cdot}{\cdot}_{\spH},\norm{\cdot}_{\spH})$
and
$(\spK, 
\ip{\cdot}{\cdot}_{\spK},\norm{\cdot}_{\spK})$
be finite dimensional real Hilbert spaces.
For $S\subset \spH$,
$\ri (S)$ denotes the relative interior of $S$.
$\setLO{\spH}{\spK}$
denotes the set of all linear operators
from 
$\spH$
to
$\spK$.
For $L \in \setLO{\spH}{\spK}$,
$\norm{L}_\mathrm{op}$ denotes the operator
norm of $L$ (i.e.,
$\norm{L}_\mathrm{op}
:= \sup_{x\in\spH, \norm{x}_\spH \leq 1}
\norm{Lx}_\spK$)
and
$L^*\in\setLO{\spK}{\spH}$ the adjoint operator
of $L$ 
(i.e., $(\forall x \in \spH) (\forall y \in \spK)\  
\ip{L x}{y}_\spK = \ip{x}{L^* y}_\spH$).
Let us denote the identity operator
of general Hilbert spaces by $\Id$ 
and the zero operators from $\spH$ to $\spK$ 
and from $\spH$ to $\spH$
by $\zeroMatrix_{\setLO{\spH}{\spK}}$
and $\zeroMatrix_{\spH}$, respectively.
For $L\in \setLO{\spH}{\spK}$, 
$\ran(L):= \{Lx\in\spK| x\in\spH \}$
and
$\nullsp(L):=\{x\in\spH | Lx=0_\spK\}$,
where
$0_\spK$ stands for the zero vector in $\spK$.
We express the positive 
definiteness and the positive semi-definiteness 
of a self-adjoint operator $L\in \setLO{\spH}{\spH}$
respectively
by $L\succ \zeroMatrix_\spH$ and 
$L\succeq \zeroMatrix_\spH$.
Any $L\succ \zeroMatrix_\spH$
defines
a new real Hilbert space
$(\spH, \ip{\cdot}{\cdot}_L, \norm{\cdot}_L)$
equipped with
an inner product $\ip{\cdot}{\cdot}_L:
\spH\times\spH\to \setR:(x,y)\mapsto\ip{x}{Ly}_\spH$
and its induced norm 
$\norm{\cdot}_L: \spH\to\setNNR: 
x\mapsto \sqrt{\ip{x}{x}_L}$.
For a matrix $A$,
$A^\transT$ denotes the transpose of $A$ and
$A^\dagger$ the Moore-Penrose pseudoinverse
of $A$.
We use $\identityMatrix_n\in \setR^{n\times n}$ to denote
the identity matrix.
We also use $\zeroMatrix_{m\times n}\in \setR^{m\times n}$
for the zero matrix.

\section{cLiGME model and its proximal splitting algorithm}
\label{sec: cLiGME model}
We start with simple reformulations 
of \eqref{eq: LiGME with convex constraints}.
\subsection{cLiGME model}
\begin{problem}[cLiGME model]
  \label{def: LiGME with C}
  Let
  $(\mathcal{X},
  \ip{\cdot}{\cdot}_{\mathcal{X}},\norm{\cdot}_{\mathcal{X}})$,
  $(\mathcal{Y}, 
  \ip{\cdot}{\cdot}_{\mathcal{Y}},\norm{\cdot}_{\mathcal{Y}})$,
  $(\mathcal{Z},
  \ip{\cdot}{\cdot}_{\mathcal{Z}},\norm{\cdot}_{\mathcal{Z}})$,
  $(\mathfrak{Z},
  \ip{\cdot}{\cdot}_{\mathfrak{Z}},\norm{\cdot}_{\mathfrak{Z}})$
  and 
  $\left( \sptildeZ,
  \ip{\cdot}{\cdot}_{\sptildeZ},
  \norm{\cdot}_{\sptildeZ}
  \right)$
  be finite dimensional real Hilbert spaces
  and 
  $ \mathbf{C}(\subset\spfrakZ)$ be a nonempty closed convex set.
  Let
  $\Psi\in \Gamma_0(\mathcal{Z})$ be coercive with
  $\dom\Psi = \spZ$. 
  Let
  $(A,B,\opL,\opC,\mu) 
  \in \mathcal{B}(\spX,\spY)
  \times \setLO{\spZ}{\sptildeZ}
  \times\mathcal{B}(\spX,\spZ)
  \times \mathcal{B}(\spX,\spfrakZ)
  \times \setPR$
  satisfy $\mathbf{C} \cap \ran \opC \neq \emptyset$.
  Then, we consider a convexly constrained LiGME (cLiGME) model:
  \begin{align}
    \label{eq: LiGME with C}
    \minimize_{\opC x\in \mathbf{C}} 
    \frac{1}{2}\norm{y -Ax}^2_{\mathcal{Y}} +
    \mu \Psi_{B} \circ \opL(x).
  \end{align}
\end{problem}

\begin{remark}[The model 
  \eqref{eq: LiGME with convex constraints}
  is a specialization of the model \eqref{eq: LiGME with C}]
  \label{remark: multiple constraints}
  By using 
  $(\spfrakZ_i,\opC_i,C_i)\ (i\in I)$
  in \eqref{eq: LiGME with convex constraints},
  define a new real Hilbert space\footnote{
  A new Hilbert space 
  $\spfrakZ := \bigtimes_{i\in I} \spfrakZ_i$
  is equipped with
  the addition
  $\spfrakZ\times\spfrakZ \to\spfrakZ:
  (x,y)=((x_i)_{i\in I}, (y_i)_{i\in I})\mapsto
  (x_i +y_i)_{i\in I}$,
  the scalar product
  $\setR\times \spfrakZ\to\spfrakZ: 
  (\alpha, (x_i)_{i\in I}) \mapsto (\alpha x_i)_{i\in I}
  $ and
  the inner product
  $\spfrakZ\times\spfrakZ \to\setR:
  (x,y) \mapsto \sum_{i\in I}\ip{x_i}{y_i}_{\spfrakZ_i}$.
  }
  $\spfrakZ := \bigtimes_{i\in I} \spfrakZ_i$,
  a new linear operator
  $\opC:\spX\to\spfrakZ: x\mapsto (\opC_i x)_{i\in I}$ and
  a new closed convex set
  $\mathbf{C} := \bigtimes_{i\in I}C_i\subset \spfrakZ$.
  Then,
  we have the equivalence
  \begin{align*}
    \opC x\in \mathbf{C}
    \iff
    (\forall i \in I)\ \opC_i x\in C_i.
  \end{align*}
\end{remark}
We use the following 
constraint-free reformulation of \eqref{eq: LiGME with C}.
\begin{lemma}[A constraint-free reformulation of cLiGME model]\label{def: cLiGME}
  In Problem \ref{def: LiGME with C},
  define new Hilbert spaces 
  $\spZ_c:= \spZ\times \spfrakZ$ and 
  $\sptildeZ_c:= \sptildeZ\times \spfrakZ$,
  new linear operators
  $\opLC:\spX\to\spZ_c:
  x\mapsto(\opL(x),\opC(x))$ and
  $B_c:=B\oplus \zeroMatrix_{\spfrakZ}:
  \spZ_c\to\sptildeZ_c:(z_1,z_2)
  \mapsto (B z_1, 0_{\spfrakZ})$,
  and a new convex function
  $\bfPsi:=\Psi\oplus\iota_{\mathbf{C}}\in \Gamma_0(\spZ_c)$
  by $\bfPsi:
  \spZ\times \spfrakZ \to (-\infty,\infty]:
  (z_1,z_2) \mapsto \Psi(z_1)+\iota_{\mathbf{C}}(z_2)$.
  Then, the cLiGME model \eqref{eq: LiGME with C}
  can be formulated equivalently as
  \begin{align}
    \label{cLiGME}
    \hspace*{-2.5mm}
    \minimize_{x\in \mathcal{X}} J_{\bfPsi_{B_c}\circ \opLC}(x):=
    \frac{1}{2}\norm{y -Ax}^2_{\mathcal{Y}} +
    \mu \bfPsi_{B_c} \circ \opLC(x),
  \end{align}
  where
  \begin{align}
  \bfPsi_{B_c}(\cdot) &:= \bfPsi(\cdot) -
  \min_{v\in \sptildeZ_c} \left[
    \bfPsi(v) 
    + \frac{1}{2}\norm{B_c (\cdot - v)}_{\spZ_c}^2
  \right]\\
  &=\Psi_B(\cdot)+\iota_\mathbf{C}(\cdot).
\end{align}
\end{lemma}

The overall convexity condition 
for the cLiGME model is given below.
\begin{proposition}[Overall convexity condition for the cLiGME model]
  \label{prop: overall convexity}
  Let $(A,\opL,\mu) 
    \in \mathcal{B}(\mathcal{X},\mathcal{Y})
    \times\mathcal{B}(\mathcal{X},\mathcal{Z})\times \mathbb{R}_{++}$.
    For three
    conditions
    $(C_1)\ A^*A-\mu\opL^*
    B^*B \opL
    \succeq \zeroMatrix_\mathcal{X}$,
    $(C_2)\ J_{\bfPsi_{B_c}\circ \opLC}\in\Gamma_0(\mathcal{X})
    \mbox{ for any }y\in \mathcal{Y}$,
    $(C_3)\ J_{\bfPsi_{B_c}\circ \opLC}^{(0)}:=
    \frac{1}{2}\norm{A\cdot}^2_\mathcal{Y}+
    \mu\bfPsi_{B_c}\circ\opLC\in\Gamma_0(\mathcal{X})$
    , the relation $(C_1)\implies (C_2)\iff (C_3)$ holds.
\end{proposition}

\subsection{A proximal splitting algorithm for cLiGME model}
Our target is the following convex optimization problem.
\begin{problem}\label{prob: cLiGME}
  Suppose $\Psi\in \Gamma_0$ satisfies
  even symmetry
  $\Psi \circ(-\Id)=\Psi$
  and is proximable (i.e., $\prox_{\gamma\Psi}$ is available as a
  computable operator for every
  $\gamma\in\mathbb{R}_{++}$)
  and assume $\mathbf{C}\subset \spfrakZ$
  is a closed convex
  subset onto which the metric projection\footnote{
    The metric projection onto a closed convex subset $\mathbf{C}\subset\spfrakZ$
    is defined by 
    \begin{align*}
      P_{\mathbf{C}}:\spfrakZ\to \mathbf{C}:z\mapsto \argmin_{v\in \mathbf{C}}\norm{z-v}_{\spfrakZ}.
    \end{align*}
    By this definition, we have 
    $\prox_{\iota_{\mathbf{C}}} = P_{\mathbf{C}}$.
  } 
  $P_{\mathbf{C}}(x):\spfrakZ\mapsto \mathbf{C}$ is computable.
  Then, for
  $\left(A,\opL,\opC,B_c,y,\mu \right)
  \in \mathcal{B}(\mathcal{X},\mathcal{Y})
  \times \mathcal{B}(\mathcal{X},\mathcal{Z})
  \times \setLO{\spX}{\spfrakZ}
  \times \mathcal{B}(\spZ_c,\sptildeZ_c)
  \times \mathcal{Y}
  \times \mathbb{R}_{++}$ satisfying
  $A^*A-\mu\opL B^*B\opL
  \succeq \zeroMatrix_\mathcal{X}$
  and $0_{\spfrakZ}\in \ri(\mathbf{C} - \ran \opC)$,
  \begin{align*}
    \mbox{\upshape find }\xstar\in \mathcal{S}
    :=\argmin_{x\in \mathcal{X}}J_{\bfPsi_{B_c}\circ \opLC}(x)
    = \argmin_{\opC x\in \mathbf{C}} 
    J_{\Psi_{B}\circ \opL} (x).
  \end{align*}
\end{problem}
The next theorem
presents an iterative algorithm of 
guaranteed convergence
to a global minimizer
of Problem \ref{prob: cLiGME}.
\begin{theorem}[Averaged nonexpansive operator $\TCC$
  and its fixed point approximation]
  \label{thm: TCC}
  In Problem \ref{prob: cLiGME},
  let $\left(\mathcal{H}:=\mathcal{X}\times \spZ \times\spZ_c,
  \ip{\cdot}{\cdot}_\mathcal{H},\norm{\cdot}_\mathcal{H}\right)$
  be a real Hilbert space, 
  Define
  $\TCC:\mathcal{H}\to\mathcal{H}:
  (x,v,w)\mapsto(\xi,\zeta,\eta) 
  \mbox{ with }(\sigma,\tau)\in \mathbb{R}_{++}\times \mathbb{R}_{++}$
  by
  \begin{align*}
    &\xi :=\left[\Id - 
    \frac{1}{\sigma}\left(A^*A 
    - \mu\opL^* B^*B\opL\right)
    \right]x
    -\frac{\mu}{\sigma}\opL^* B^*B v\\
    &\quad\quad -\frac{\mu}{\sigma} \opLC^*w+A^*y\\
    &\zeta := \prox_{\frac{\mu}{\tau}\Psi}\left[
      \frac{2\mu}{\tau} B^*B \opL\xi
      - \frac{\mu}{\tau}B^*B\opL x
      +\left(\Id - \frac{\mu}{\tau}B^*B\right)v
    \right]\\
    &\eta := 
    (\Id - \prox_{\Psi \oplus \iota_{\mathbf{C}}})
    \left(2\opLC\xi-\opLC x+w\right),
  \end{align*}
  where 
  $\prox_{\Psi\oplus \iota_{\mathbf{C}}}(w_1,w_2) 
  =(\prox_{\Psi}(w_1),P_{\mathbf{C}}(w_2))$.
  Then,
  \begin{enumerate}[\upshape (a)]
    \item The solution set $\mathcal{S}$ of Problem 
    \ref{prob: cLiGME} can be expressed as
    \begin{align*}
      \mathcal{S}&=\Xi\left( \fix\left( \TCC\right) \right)\\
      &:=\left\{ \Xi\left(x, v, w \right)
      | (x,v,w)\in \fix(\TCC)\right\},
    \end{align*}
    where $\Xi:\mathcal{H}\to\mathcal{X}:
    \left(x, v, w \right)\mapsto x$
    and $\fix(\TCC):= \{
      (x,v,w)\in \spH | \left(x, v, w \right) = \TCC 
      \left(x, v, w \right)
      \}$.
    \item
    Choose 
    $(\sigma,\tau,\kappa)
      \in \mathbb{R}_{++}\times \mathbb{R}_{++} \times(1,\infty)$
    satisfying\footnote{
      For example, 
    choose any $\kappa > 1$. Compute $(\sigma, \tau)$ by
    $\sigma:= \norm{\frac{\kappa}{2}A^*A
    +\mu\opLC^*\opLC}_{\mathrm{op}}+(\kappa - 1)$
    and 
    $\tau := \left(\frac{\kappa}{2} + \frac{2}{\kappa}\right)\mu\norm{B}_{\mathrm{op}}^2+(\kappa-1)$
    Then, $(\sigma,\tau,\kappa)$ satisfies
    the condition \eqref{eq: cond for definiteness}.
    }
      \begin{align}
        \begin{cases}
          &\sigma \Id -\frac{\kappa}{2}A^*A-\mu \opLC^*\opLC
          \succ \zeroMatrix_\spX\\
          &\tau \geq\left( \frac{\kappa}{2}+ \frac{2}{\kappa}
          \right)\mu \norm{B}_{\mathrm{op}}^2.
        \end{cases}\label{eq: cond for definiteness}
      \end{align}
      Then,
      \begin{align}
        \opP :=\begin{bmatrix}
          \sigma \Id & -\mu\opL^*B^*B & -\mu \opLC^*\\
          -\mu B^*B\opL& 
          \tau\Id & \zeroMatrix_{\setLO{\spZ_c}{\spZ}}\\
          -\mu \opLC & \zeroMatrix_{\setLO{\spZ}{\spZ_c}} & \mu \Id
        \end{bmatrix}\succ\zeroMatrix_\spH \label{def: operatorB}
      \end{align}
      holds and
      $\TCC$ is $\frac{\kappa}{2\kappa-1}$-averaged nonexpansive 
      in $(\spH,\ip{\cdot}{\cdot}_\opP,\norm{\cdot}_\opP)$
      which means that there exists
      $\hat{T}:\spH\to\spH$ satisfying 
      \begin{align*}
        (\forall \mathbf{x},\mathbf{y} \in \spH)
        \norm{\hat{T}\mathbf{x}-\hat{T}\mathbf{y}}_{\opP}
        \leq
        \norm{\mathbf{x}-\mathbf{y}}_{\opP}
      \end{align*}
      and $\TCC = \left(1-\frac{\kappa}{2\kappa - 1}\right)
      \Id + \frac{\kappa}{2\kappa - 1}
      \hat{T}$.

      \item
      
      Assume that
      $(\sigma,\tau,\kappa)
      \in \mathbb{R}_{++}\times \mathbb{R}_{++} \times(1,\infty)$
      satisfies \textup{(\ref{eq: cond for definiteness})}.
      Then, for any initial point $(x_0,v_0,w_0)\in \spH$,
      the sequence 
      $(x_k,v_k,w_k)_{k\in\setN}$
      generated by
      \begin{align*}
        (\forall k\in\setN)\ (x_{k+1},v_{k+1},w_{k+1})
        = \TCC(x_k,v_k,w_k)
      \end{align*}     
      converges to a point 
      $(x^\diamond,v^\diamond,w^\diamond)\in \fix(\TCC)$ and
      \begin{align*}
        \lim_{k\to\infty}x_k=x^\diamond\in \setS.
      \end{align*}
  \end{enumerate}
\end{theorem}
Algorithm \ref{alg for prob1} below is made based
on Theorem \ref{thm: TCC}.
\begin{algorithm}[H]
  \caption{for Problem (\ref{prob: cLiGME})}
  \label{alg for prob1}
  \begin{algorithmic}
  {\small
    \State Choose $(x_0,v_0,w_0)\in \spH$.
    \State Let $(\sigma,\tau,\kappa)
    \in \mathbb{R}_{++}\times \mathbb{R}_{++} \times(1,\infty)$ 
    satisfying (\ref{eq: cond for definiteness}).
    \State Define $\opP$ as (\ref{def: operatorB}).
    \State $k\leftarrow 0$.
    \Do
    \State $x_{k+1} \leftarrow \left[\Id - 
      \frac{1}{\sigma}\left(A^*A 
      - \mu\opL^*B^*B\opL\right)
      \right]x_k
      -\frac{\mu}{\sigma}\opL^*B^*B v_k$\\
      $\quad\quad\quad\quad\quad \  -\frac{\mu}{\sigma} \opL_c^*w_k+A^*y$
    \State $ v_{k+1} \leftarrow \prox_{\frac{\mu}{\tau}\Psi}[
      \frac{2\mu}{\tau}B^*B \opL x_{k+1}
      - \frac{\mu}{\tau}B^*B\opL x_k$\\
      $\quad\quad\quad\quad\quad\quad\quad\quad \quad
      \
      +(\Id - \frac{\mu}{\tau}B^*B)v_k
      ]$
    \State $w_{k+1} \leftarrow 
    (\Id - \prox_{\Psi \oplus \iota_{\mathbf{C}}})
    \left(2\opL_c x_{k+1}-
      \opL_c x_k+w_k\right)$
    \State $k\leftarrow k+1$ 
    \doWhile{$\norm{(x_k,v_k,w_k)
    - (x_{k-1},v_{k-1},w_{k-1})}_{\opP}$
    is not sufficiently small}\\
    \Return $x_k$
  }
  \end{algorithmic}
\end{algorithm}

\subsection{How to choose $B$ to achieve the overall convexity}
A matrix $B$ achieving the overall convexity of
$J_{\bfPsi_{B_c} \circ \opL_c}$
can be designed in exactly same way as 
\cite[Proposition 2]{LiGME} for \eqref{cLiGME}.
\begin{proposition}[A design of $B$ to ensure the 
  overall convexity condition in Proposition 
  \ref{prop: overall convexity}]
  \label{prop: how to choose B}
  In Problem \ref{def: LiGME with C}, let 
  $(\spX,\spY,\spZ) = (\setR^n, \setR^m,\setR^l)$,
  $(A,\opL,\mu) 
  \in\setR^{m\times n}\times \setR^{l\times n}
  \times \setPR$
  and $\rank(\opL) = l$. Choose a nonsingular
  $\tildeL\in \setR^{n \times n}$ satisfying
  $\begin{bmatrix}
    \zeroMatrix_{l \times (n-l)} & \identityMatrix_l
  \end{bmatrix}\tildeL = \opL$.
  Then,
  \begin{align}
    B = B_\theta := \sqrt{\theta / \mu}\Lambda^{1/2} U^{\transT}
    \in \setR^{l\times l}, \theta \in [0,1]
  \end{align}
  ensures $J_{\bfPsi_{B_c}\circ\opL_c}\in \Gamma_0(\setR^n)$,
  where 
  \begin{align}
    \begin{bmatrix}
      \tilde{A_1} & \tilde{A_2}
    \end{bmatrix} = A(\tildeL)^{-1} \label{eq:def of A1 A2}
  \end{align}
  and $U\Lambda U^\transT :=\tilde{A_2}^\transT \tilde{A_2}
  -\tilde{A_2}^\transT \tilde{A_1} \left(\tilde{A_1}^\transT
  \tilde{A_1}\right)^\dagger \tilde{A_1}^\transT \tilde{A_2}$
  is an eigendecomposition.
\end{proposition}

In the following proposition,
(a) shows that the cLiGME model \eqref{eq: LiGME with C}
admits multiple penalties 
and (b) shows how to design $B_\theta$ for multiple penalties.

\begin{proposition}[Multiple penalties for the cLiGME model]
  \quad \label{prop: multiple penalties}
  \begin{enumerate}[\upshape (a)]
    \item \label{prop: reduce multiple penalties}
    \textup{(The sum of multiple LiGME penalties
    can be modeled as a single cLiGME penalty on product space).}
      Let $\spZ_i$ and $\sptildeZ_i$
      be finite dimensional real Hilbert spaces 
      $(i = 1,2,3,\cdots \mathcal{M})$,
      $\Psi^{\lrangle{i}} \in \Gamma_0(\spZ_i)$
      be coercive
      with $\dom(\Psi^{\lrangle{i}}) = \spZ_i$,
      $\mu_i\in \setPR$, 
      $B^{\lrangle{i}}\in \setLO{\spZ_i}{\sptildeZ_i}$,
      and $\opL_i\in \setLO{\spZ_i}{\sptildeZ_i}$.
      Define a new real Hilbert spaces
      $\spZ := \bigtimes_{i = 1}^{\mathcal{M}}\spZ_i$
      and
      $\sptildeZ := \bigtimes_{i = 1}^{\mathcal{M}}\sptildeZ_i$,
      a new linear operator 
      $B: \spZ\to \sptildeZ:(z_1,\cdots,z_{\mathcal{M}})
      \mapsto (\sqrt{\mu_1}B^{\lrangle{1}}z_1,\cdots,
      \sqrt{\mu_{\mathcal{M}}}B^{\lrangle{\mathcal{M}}}
      z_{\mathcal{M}})$
      and
      $\opL:\spX\to\spZ: x\mapsto(
        \opL_i x_i
      )_{1\leq i \leq\mathcal{M}}$,
      and a new function
      $\Psi :=\bigoplus_{i = 1}^{\mathcal{M}}\mu_i 
      \Psi^{\lrangle{i}}$.
      Then, we have
      \begin{align}
        \label{eq: multiple penalties}
        \bfPsi_{B_c}\circ\opLC= \left(\sum_{i=1}^{\mathcal{M}}
        \mu_i\left(
          \Psi^{\lrangle{i}}
        \right)_{B^{\lrangle{i}}}\circ\opL_i
        \right)
        + (\iota_\mathbf{C}\circ \opC).
      \end{align}

    \item \label{prop: multiple B}
    \textup{(A design of $B^{\lrangle{i}}$ to achieve 
    overall convexity in Prop. 
    \ref{prop: overall convexity}).}
    Let $(\spX,\spY,\spZ_i) = (\setR^n,\setR^m,\setR^{l_i})$,
    $(A,\opL_i,\mu)\in \setLO{\setR^n}{\setR^m}\times
    \setLO{\setR^{l_i}}{\setR^n}\times\setPR$
    and 
    $\rank(\opL_i)= l_i$ ($i = 1,2,3,\cdots \mathcal{M}$).
    For each $i = 1,2,\cdots \mathcal{M}$,
    choose nonsingular $\tildeL_i\in \setR^{n\times n}$
    satisfying 
    $\begin{bmatrix}
      \zeroMatrix_{l\times (n-l)} & \identityMatrix_l
    \end{bmatrix}\tildeL_i = \opL$ and $\omega_i\in\setPR$
    satisfying
    $\sum_{i=1}^\mathcal{M} \omega_i=1$.
    For each $i=1,2,\cdots, \mathcal{M}$,
    apply Proposition \ref{prop: how to choose B}
    to $\left(\sqrt{\frac{\omega_i}{\mu}}A, \opL_i, \mu_i\right)$
    to obtain $B_{\theta_i}^{\lrangle{i}}\in \setR^{l_i\times l_i}$
    satisfying
    $\left(\sqrt{\frac{\omega_i}{\mu}}
    A\right)^\transT\left(
      \sqrt{\frac{\omega_i}{\mu}}A
    \right)-\mu_i \opL_i^\transT
    \left(B_{\theta_i}^{\lrangle{i}}\right)^\transT
    B_{\theta_i}^{\lrangle{i}}\opL_i\succeq \zeroMatrix_{n\times n}$. 
    Then,
    for 
    $\theta = (\theta_1,\theta_2,$\\
    $\cdots,\theta_{\mathcal{M}})$, $B_\theta:
    \bigtimes_{i=1}^{\mathcal{M}}\setR^{l_i}
    \to\bigtimes_{i=1}^{\mathcal{M}}\setR^{l_i}:
    (z_1,\cdots,z_{\mathcal{M}})
    \mapsto
    (\sqrt{\mu_1} B_{\theta_1}^{\lrangle{1}}z_1,
    \cdots,
    \sqrt{\mu_{\mathcal{M}}}
    B_{\theta_{\mathcal{M}}}^{\lrangle{\mathcal{M}}}z_{\mathcal{M}}
    )$
    ensures $J_{\bfPsi_{B_c}\circ\opLC}\in \Gamma_0(\setR^n)$
    (see Lemma \ref{def: cLiGME}).
  \end{enumerate}
\end{proposition}

\section{Numerical experiments}
\label{sec: numerical experiment}
\begin{figure}[h]
  \centering
  \begin{minipage}[b]{0.24\linewidth}
    \centering
    \includegraphics[keepaspectratio, scale=0.1]
    {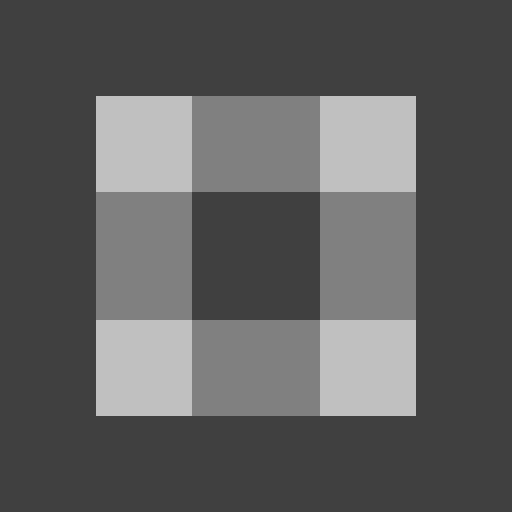}
    \subcaption{}
  \end{minipage}
  \begin{minipage}[b]{0.24\linewidth}
    \centering
    \includegraphics[keepaspectratio, scale=0.1]
    {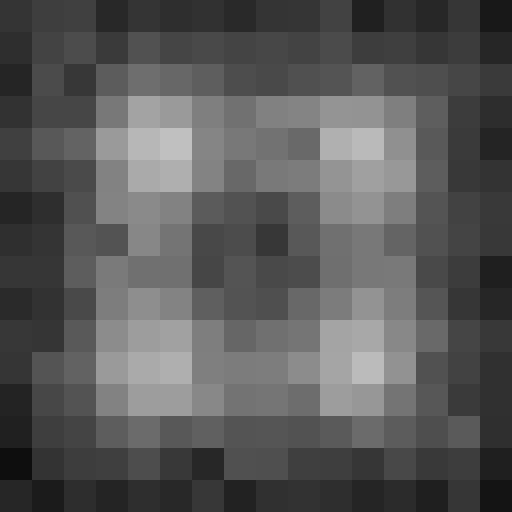}
    \subcaption{}
  \end{minipage}\\
  \begin{minipage}[b]{0.242\linewidth}
    \centering
    \includegraphics[keepaspectratio, scale=0.1]
    {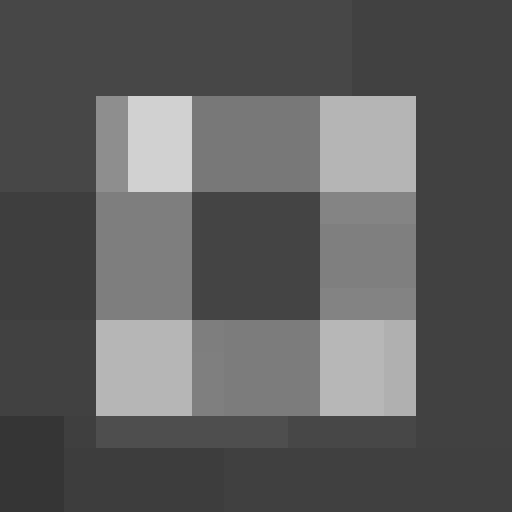}
    \subcaption{TV $\mathbf{C}_{\clubsuit}$}
  \end{minipage}
  \begin{minipage}[b]{0.24\linewidth}
    \centering
    \includegraphics[keepaspectratio, scale=0.1]
    {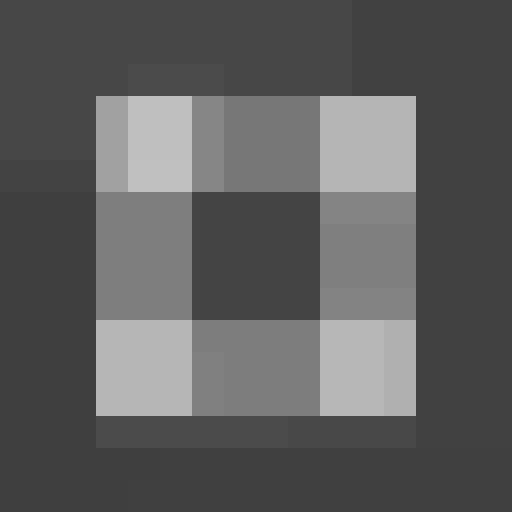}
    \subcaption{TV $\mathbf{C}_{\diamondsuit}$}
  \end{minipage}
  \begin{minipage}[b]{0.24\linewidth}
    \centering
    \includegraphics[keepaspectratio, scale=0.1]
    {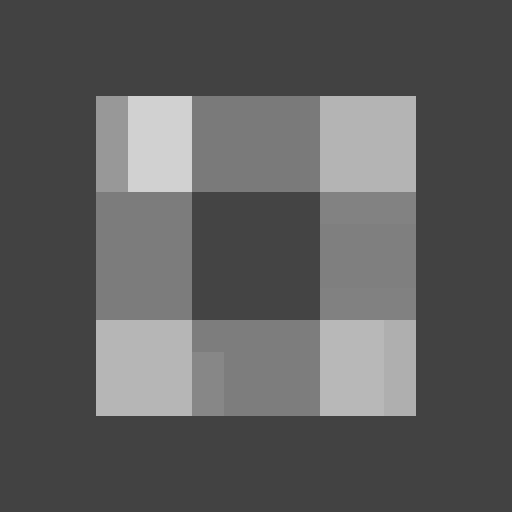}
    \subcaption{TV $\mathbf{C}_{\heartsuit}$}
  \end{minipage}
  \begin{minipage}[b]{0.24\linewidth}
    \centering
    \includegraphics[keepaspectratio, scale=0.1]
    {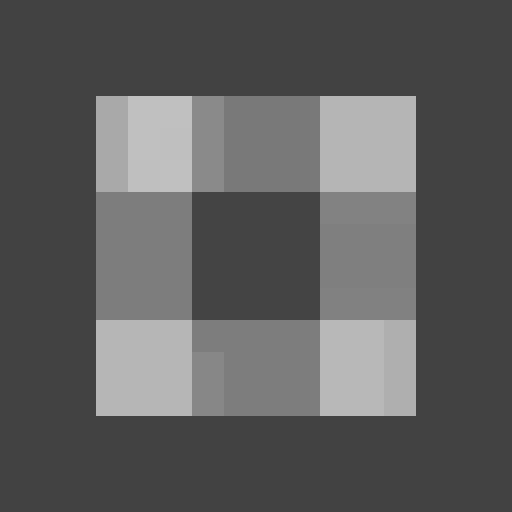}
    \subcaption{TV $\mathbf{C}_\spadesuit$}
  \end{minipage}\\
  \begin{minipage}[b]{0.24\linewidth}
    \centering
    \includegraphics[keepaspectratio, scale=0.1]
    {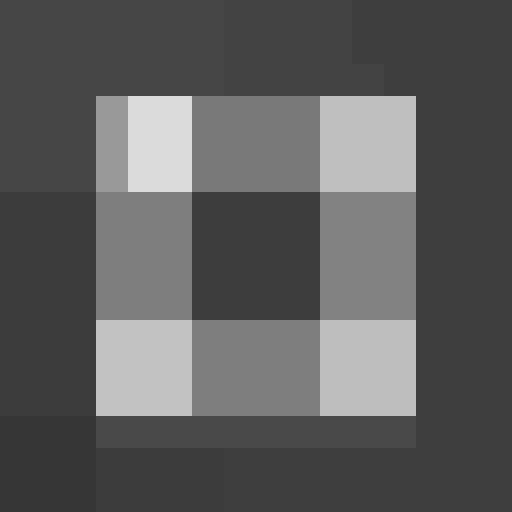}
    \subcaption{cLiGME $\mathbf{C}_{\clubsuit}$}
  \end{minipage}
  \begin{minipage}[b]{0.24\linewidth}
    \centering
    \includegraphics[keepaspectratio, scale=0.1]
    {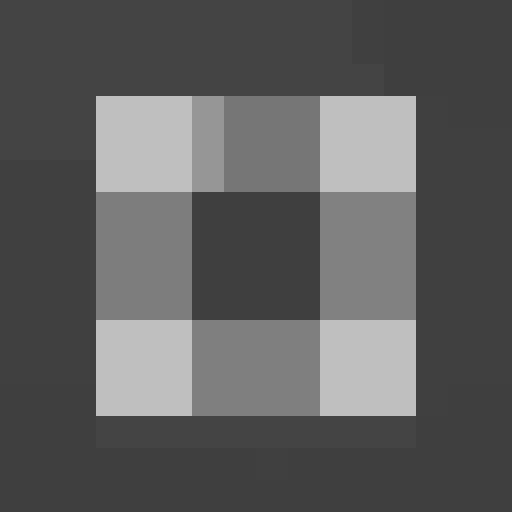}
    \subcaption{cLiGME $\mathbf{C}_{\diamondsuit}$}
  \end{minipage}
  \begin{minipage}[b]{0.24\linewidth}
    \centering
    \includegraphics[keepaspectratio, scale=0.1]
    {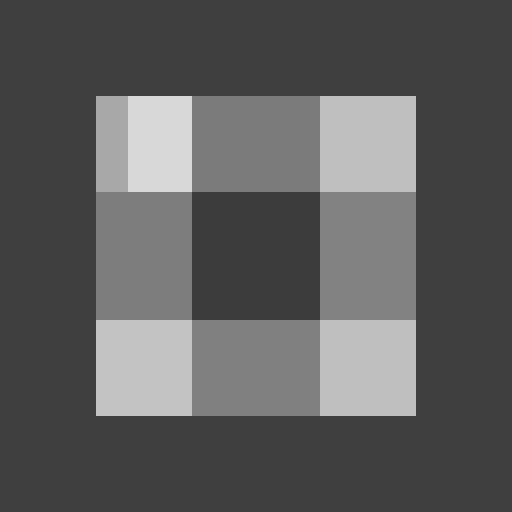}
    \subcaption{cLiGME $\mathbf{C}_{\heartsuit}$}
  \end{minipage}
  \begin{minipage}[b]{0.24\linewidth}
    \centering
    \includegraphics[keepaspectratio, scale=0.1]
    {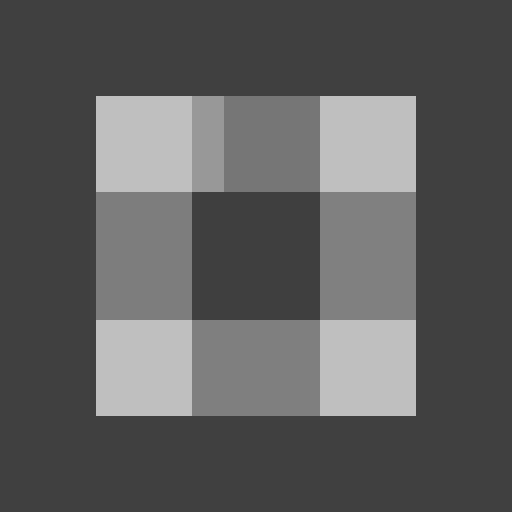}
    \subcaption{cLiGME $\mathbf{C}_{\spadesuit}$}
  \end{minipage}
  \caption{(a) original piecewise constant image of
  pixel values in $\{0.25,0.5,0.75\}$,
  (b) observed noisy blurred image,
  (c-j) recovered images.
  Each recovered image is obtained by 
  Algorithm \ref{alg for prob1}
  at $5000$th iteration 
  and shown with the pixel value range 
  $[0 (\mbox{black}),1 (\mbox{white})]$.
  }
  \label{pic: images of X}
\end{figure}

We conducted numerical experiments based on 
a scenario of image restoration for
piecewise constant $N$-by-$N$ image, which is 
the same as \cite[Sec. 4.2]{LiGME}
for space saving purpose but with considering
additionally multiple convex constraints.
Set $(\spX,\spY,\spZ_1,\spZ_2)=
\left(\setR^{N^2}, \setR^{N^2},\setR^{N(N-1)},
\setR^{N(N-1)}\right)$, $N = 16$,
$\Psi^{\lrangle{1}}= \Psi^{\lrangle{2}}= \norm{\cdot}_1$
, $\mu_1 = \mu_2 = 1$ and 
$\opL =\bar{D}:= [\dmath^\transT, \dmatv^\transT]^\transT$,
where $\dmatv\in \setR^{N(N-1)\times N^2}$ is the
vertical difference operator 
and $\dmath \in \setR^{N(N-1)\times N^2}$ 
is the horizontal difference operator used
in \cite[(43)]{LiGME}, respectively 
(see Proposition \ref{prop: multiple penalties}(\ref{prop: reduce multiple penalties})).
The blur matrix $A\in \setR^{N^2\times N^2}$ is also defined as in \cite[(44) and (45)]{LiGME}.
The observation vector $y\in \setR^{N^2}$
is assumed to satisfy the linear regression model:
$y = A\xstar + \varepsilon$,
where $\xstar\in \setR^{N^2}$ is 
the vectorization\footnote{The vectorization of a matrix is 
the mapping:
\begin{align*}
  \vecrize: \setR^{m\times n}\to\setR^{mn}
  :A \mapsto \left[
    a_1^\transT,a_2^\transT\cdots a_n^\transT
  \right]^\transT,
\end{align*}
where $a_i$
is the $i$-th column vector of $A$
for each $i\in \{1,2,\cdots ,n \}.$} of
a piecewise constant image shown
in Figure \ref{pic: images of X}(a)
and $\varepsilon\in \setR^{N^2}$ is the 
additive white Gaussian noise.
As available a priori knowledge
on $\xstar$,
we use (i) every 
entry in $\xstar$ belongs to $[0.25,0.75]$, 
(ii) every entry in the background $I_{\mathrm{back}}$ is 
same 
but unknown valued (similar a priori knowledge is found, e.g.,
in
blind deconvolution
\cite{A novel blind deconvolution scheme 
for image restoration using recursive filtering}).
the signal-to-noise ratio (SNR) 
is set by 
$ 10\log_{10}\frac{\norm{x}_{\spX}^2}
{\norm{\varepsilon}_\spY^2}\ =20 [\mathrm{dB}]
$.
By letting $\spX =: \spfrakZ_1 =: \spfrakZ_2$,
$\spfrakZ :=\spfrakZ_1 \times \spfrakZ_2$,
$\opC_1 := \opC_2 := \identityMatrix_{N^2}$,
$\opC:= \left[\opC_1^\transT,
\opC_2^\transT\right]^\transT$,
\begin{align*}
  C_1 &:= \left\{x\in\setR^{N^2}=:\spfrakZ_1| 
  (i = 1,2\cdots,N^2)\  
  0.25\leq x_i \leq 0.75 
  \right\},\\
  \begin{split}
    C_2 &:= \left\{\vecrize (X)\in \setR^{N^2}=:\spfrakZ_2| 
    X\in \setR^{N\times N},\right.\\
    &\quad\quad\left.
    (\forall (i_1,j_1) \in  I_{\mathrm{back}})
    (\forall (i_2,j_2) \in  I_{\mathrm{back}})
    \ X_{i_1,j_1} = X_{i_2,j_2} 
    \vphantom{\setR^{N^2}} 
    \right\},
  \end{split}
\end{align*}
where 
$
  I_{\mathrm{back}}:=\left[\left\{
    1,2,3,14,15,16
  \right\}\times
  \left\{
    1,2,\cdots , 16
  \right\}\right]\\
  \cup
  \left[ 
  \left\{
    1,2,\cdots , 16
  \right\} \times
  \left\{
    1,2,3,14,15,16
  \right\}
  \right]
$,
we consider the following four cases: \\
($\clubsuit$) \ \
$\opC x^\star \in\mathbf{C}_{\clubsuit} := \spfrakZ :=\setR^{N^2}\times\setR^{N^2}$\\
($\diamondsuit$) \ \ 
$\opC x^\star \in\mathbf{C}_{\diamondsuit}:=C_1\times \setR^{N^2}$\\
($\heartsuit$) \ \ 
$\opC x^\star \in\mathbf{C}_{\heartsuit}:= \setR^{N^2}\times C_2$\\
($\spadesuit$) \ \ 
$\opC x^\star \in \mathbf{C}_{\spadesuit} :=C_1\times C_2$.\\
For each case ($p = \clubsuit,\diamondsuit,\heartsuit,\spadesuit$), 
we compared two models
as instances of Problem \ref{prob: cLiGME}:
\begin{enumerate}[(a)]
  \item the anisotropic TV with constraint $\opC x\in \mathbf{C}_p$:
  \begin{align}
    \label{eq: Numerical TV}
    \minimize_{\opC x\in \mathbf{C}_p} 
    J_{({\norm{\cdot}_1})_{\zeroMatrix_{\spZ_1\times
    \spZ_2}}\circ \bar{D}}(x)
    (= J_{\norm{\cdot}_1 \circ\bar{D}}(x))
  \end{align}
  \item cLiGME for 
  $(\Psi,\opL) 
  =(\norm{\cdot}_1,\bar{D})$
  with constraint $\opC x\in \mathbf{C}_p$:
  \begin{align}
    \label{eq: Numerical cLiGME}
    \minimize_{\opC x\in \mathbf{C}_p} 
    J_{({\norm{\cdot}_1})_{B_\theta}\circ \bar{D}}(x), 
  \end{align}
\end{enumerate}
where 
$\norm{\cdot}_1$ is defined on the space $\spZ_1\times\spZ_2$
and
$B_\theta$ is defined by 
Proposition \ref{prop: multiple penalties}(\ref{prop: multiple B})
with $(\theta_i, \omega_i,\tildeL_i)$ used in \cite[Sec. 4.2]{LiGME}
(Note: ($\clubsuit$) reproduces the 
experiment in \cite[Sec. 4.2]{LiGME}).

For all cases,
we applied Algorithm
\ref{alg for prob1}
with $\kappa = 1.001$,
$(\sigma,\tau)$ obtained by
footnote of Theorem \ref{thm: TCC}
and the initial point $(x_0,v_0,w_0)= (0_\spX,0_\spY,0_\spZ)$.
$\prox_{\gamma\norm{\cdot}_1}$
is computed by the soft-thresholding 
$[\prox_{\gamma \norm{\cdot}_1}]_i:\setR^{2N(N-1)}\to \setR:$
\begin{align*}
  (z_1,z_2, \cdots,z_{2N(N-1)})^\transT \mapsto
  \begin{cases}
    0 & (|z_i|\leq \gamma)\\
    (|z_i|-\gamma)\frac{z_i}{|z_i|} &(otherwise).
  \end{cases}
\end{align*}
Projections 
onto $C_1$ and $C_2$ are computed respectively by
\begin{align*}
  (P_{C_1}(x))_i &= \begin{cases}
    0.75 & (0.75 < x_i)\\
    x_i & (0.25\leq x_i \leq 0.75)\\
    0.25 & (x_i< 0.25),
  \end{cases}\\
  (P_{C_2}(x))_i &= \begin{cases}
    \frac{1}{\# I_{\mathrm{back}}}
    \sum_{\ind^{-1}(j)\in I_{\mathrm{back}}}
    x_j
    & (\ind^{-1}(i)\in I_{\mathrm{back}})\\
    x_i & (otherwise),
  \end{cases}
\end{align*}
where 
$\ind(i,j) := i+j(N-1)$.\par

\begin{figure}[h]
  \begin{minipage}[b]{0.45\linewidth}
    \centering
    \includegraphics[keepaspectratio, scale=0.11]
    {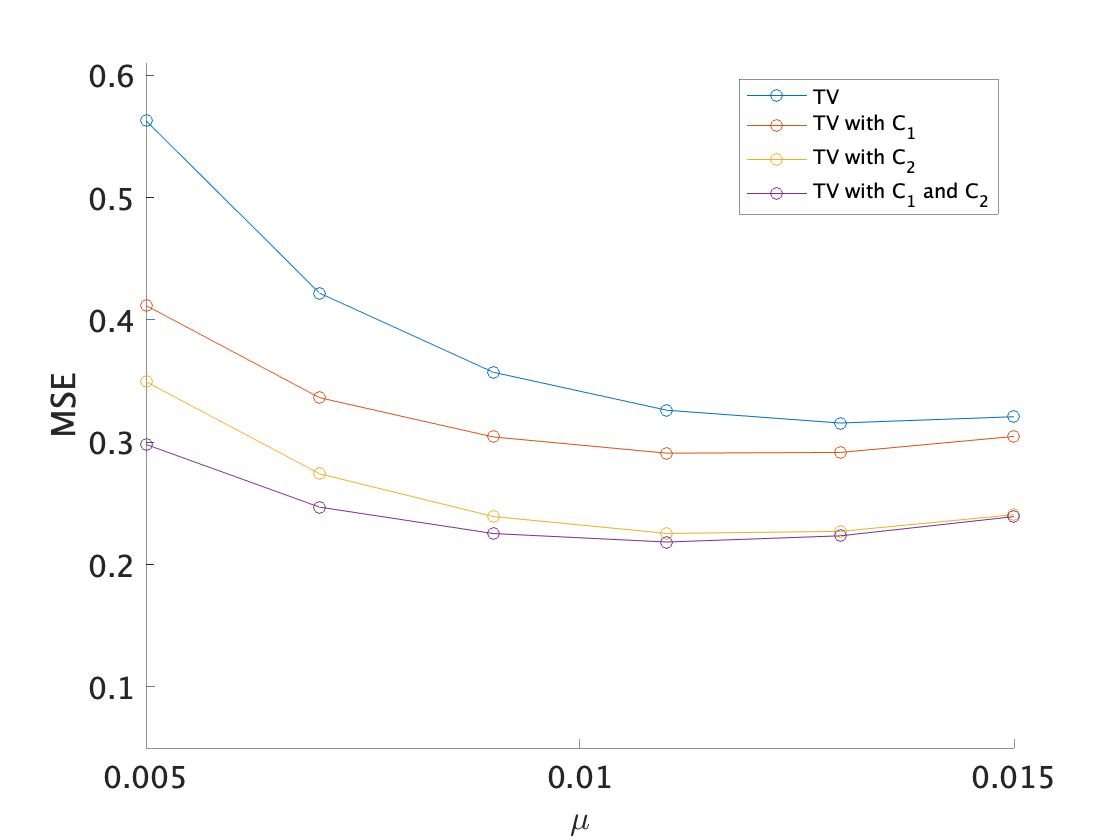}
    \subcaption{}
  \end{minipage}
  \begin{minipage}[b]{0.45\linewidth}
    \centering
    \includegraphics[keepaspectratio, scale=0.11]
    {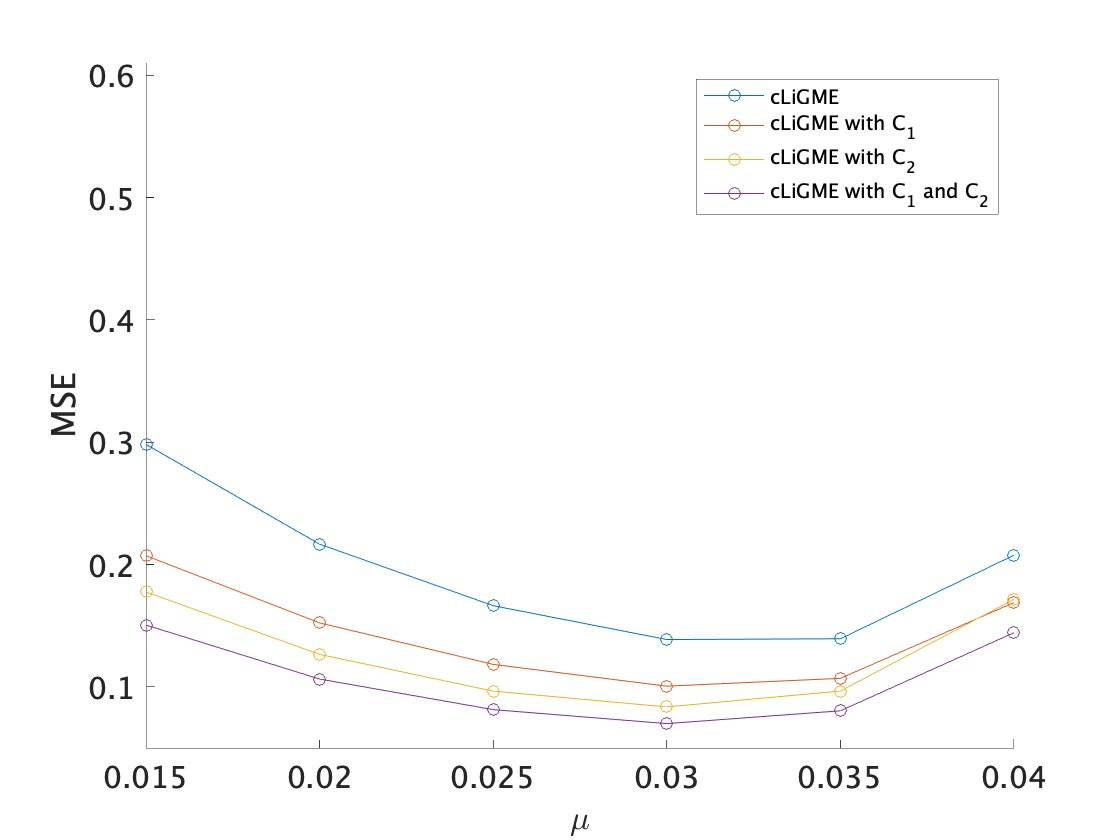}
    \subcaption{}
  \end{minipage}
  \caption{MSE versus $\mu$
  (around best $\mu$) in Problem {\ref{prob: cLiGME}}
  at $5000$th iteration for
  (a) TV model in {\eqref{eq: Numerical TV}}
  and for
  (b) cLiGME model in {\eqref{eq: Numerical cLiGME}}
  }
  \label{pic: MSE versus mu}
\end{figure}

Figure \ref{pic: MSE versus mu}
shows the dependency of recovery performance
on the parameter $\mu$ in Problem \ref{prob: cLiGME}.
The mean squared error (MSE) is defined by the
average of 
squared error (SE): 
$\norm{x_k-x^\star}_\spX^2$
over 100 independent realizations of the additive
white Gaussian noise $\varepsilon$.
We see that 
constraints are effective to improve the estimation 
and cLiGME model 
in \eqref{eq: Numerical cLiGME} can achieve better estimation
than TV model in \eqref{eq: Numerical TV}.
\par
\begin{figure}[h]
  \begin{minipage}[b]{0.45\linewidth}
    \centering
    \includegraphics[keepaspectratio, scale=0.11]
    {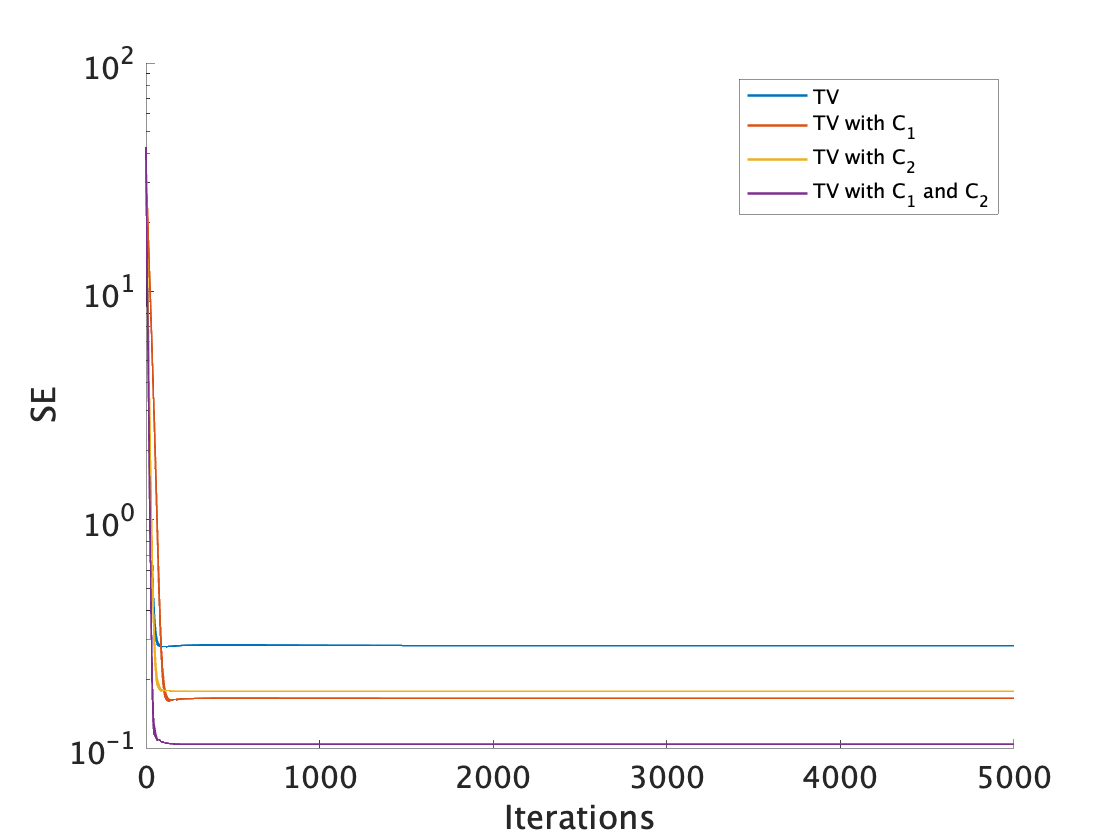}
    \subcaption{}
  \end{minipage}
  \begin{minipage}[b]{0.45\linewidth}
    \centering
    \includegraphics[keepaspectratio, scale=0.11]
    {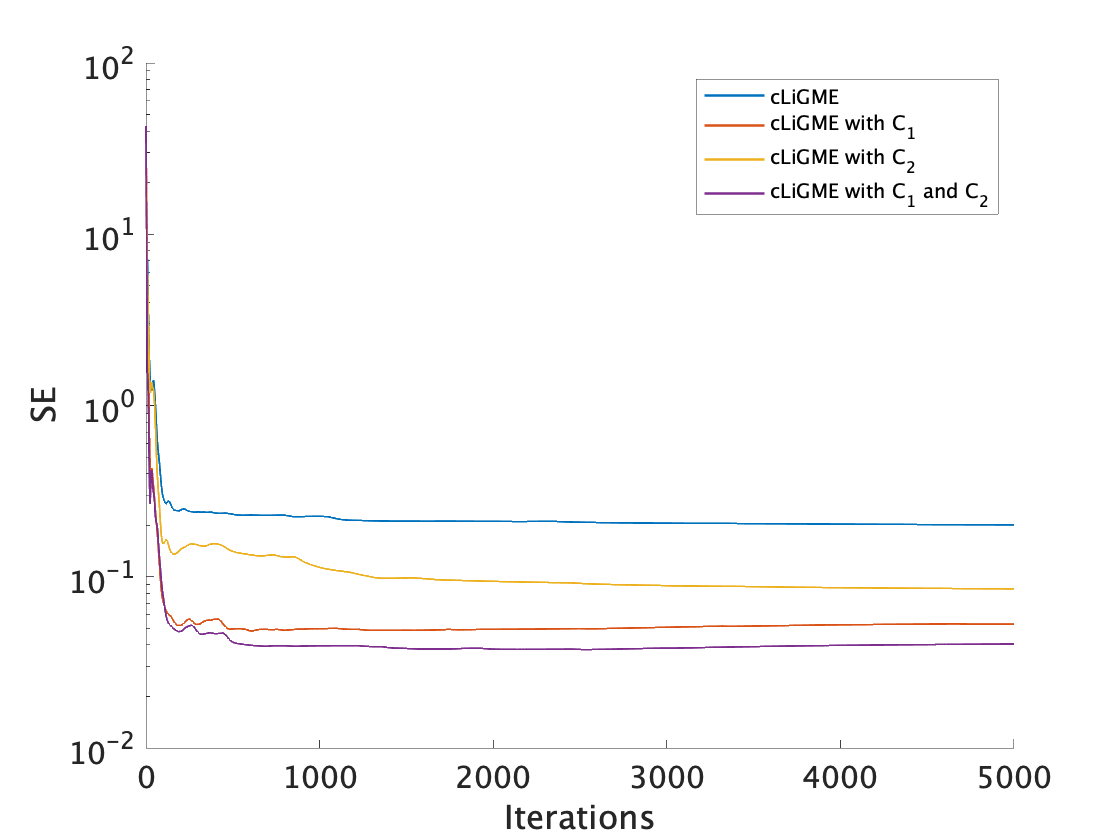}
    \subcaption{}
  \end{minipage}
  \caption{SE versus iterations
  for
  (a) TV model in {\eqref{eq: Numerical TV}}
  and for
  (b) cLiGME model in {\eqref{eq: Numerical cLiGME}}
    }\label{pic: SE versus iteration}
\end{figure}
Figure \ref{pic: SE versus iteration} shows 
convergence performances observed over $5000$ iterations.
Common weights are used as
$\mu = \mu_{\mathrm{TV}}:=0.013$ for 
\eqref{eq: Numerical TV}
and
$\mu = \mu_{\mathrm{cLiGME}}:= 0.03$ 
for \eqref{eq: Numerical cLiGME}, 
where these values
are best ones for the TV model and cLiGME model
in the case of ($\clubsuit$).
Even if we use $C_1$ and $C_2$,
the model \eqref{eq: Numerical TV} cannot achieve
the level $10^{-1}$ while
the model \eqref{eq: Numerical cLiGME} achieves
this level only with $C_1$.\par
Figure \ref{pic: images of X} shows the original image,
an observed image and recovered images by the models
\eqref{eq: Numerical TV} and \eqref{eq: Numerical cLiGME}. 

\section{Conclusion}
We proposed a convexly constrained LiGME (cLiGME)
model to broaden applicability of the LiGME model.
The cLiGME model allows to use multiple
convex constraints.
We also proposed a proximal splitting type algorithm
for cLiGME model.
Numerical experiments show 
the efficacy of the proposed model and the proposed algorithm.

\end{document}